\documentclass{elsart}
\usepackage{amscd,amsmath,amsfonts,amssymb,xy,graphicx}
\xyoption{line}
\journal{}

\newcommand{\se}[1]{\begin{equation*}\begin{split}#1\end{split}\end{equation*}}

\newcommand{\CC}{\mathbb{C}}
\newcommand{\NN}{\mathbb{N}}
\newcommand{\N}{\mathbb{N}}

\newcommand{\ZZ}{\mathbb{Z}}
\newcommand{\Z}{\mathbb{Z}}

\newcommand{\R}{\mathbb{R}}
\newcommand{\PP}{\mathbb{P}}
\newcommand{\Affine}{\mathbb{A}}
\newcommand{\vtx}[1]{*+[o][F-]{\scriptscriptstyle #1}}

\newcommand{\cR}{\mathcal{R}}

\newtheorem{theorem}{Theorem}
\newtheorem{lemma}{Lemma}

\newcommand{\prf}{\noindent\textbf{Proof.}~}
\newcommand{\eop}{\hfill $\square$\\}

\newcommand{\rep}{\texttt{rep}}
\newcommand{\trep}{\texttt{trep}}
\newcommand{\triss}{\texttt{triss}}
\newcommand{\iss}{\texttt{iss}}
\newcommand{\Flat}{\texttt{flat}}
\newcommand{\moss}{\texttt{moss}}
\newcommand{\ress}{\texttt{ress}}
\newcommand{\Null}{\texttt{Null}}
\newcommand{\Hom}{\textrm{Hom}}
\newcommand{\Image}{\textrm{Im}}
\newcommand{\GL}{\textrm{GL}}

\newcommand{\Graph}{\Gamma}
\newcommand{\BS}{\texttt{BS}}
\xyoption{all}

\begin{document}

\begin{frontmatter}
\title{First Steps towards Hyper-desingularization through Brauer-Severi Varieties}

\author{Raf Bocklandt}
\ead{rafael.bocklandt@ua.ac.be}
\author{Stijn Symens}
\ead{stijn.symens@ua.ac.be}
\author{Geert Van de Weyer\thanksref{fwo}}
\ead{geert.vandeweyer@ua.ac.be}
\thanks[fwo]{The first and the third author are Postdoctoral Fellows of the Fund for Scientific Research - Flanders (Belgium)(F.W.O. - Vlaanderen) The second author is Research Assistant of the Fund for Scientific Research - Flanders (Belgium)(F.W.O. - Vlaanderen)}
\address{Department of Mathematics and Computer Science, University of Antwerp, Middelheimlaan 1, B-2020 Antwerp, Belgium.}

\begin{abstract}
Given a Cayley-Hamilton smooth order $A$ in a central simple algebra $\Sigma$, we determine the flat locus of the Brauer-Severi fibration of the smooth order. Moreover, we give a classification of all (reduced) central singularities where the flat locus differs from the Azumaya locus and show that the fibers over the flat, non-Azumaya points near these central singularities can be described as fibered products of graphs of projection maps, thus generalizing an old result of Artin on the fibers of the Brauer-Severi fibration over a ramified point. Finally, we show these fibers are also toric quiver varieties and use this fact to compute their cohomology.
\end{abstract}
\end{frontmatter}

\section{Introduction}
Let $\CC$ be an algebraically closed field of characteristic zero and consider a function field $K$ over $\CC$ of transcendence degree $d$. Let $\Sigma$ be a central simple algebra of dimension $n^2$ with reduced trace $tr$ and let $C$ be an affine normal domain with function field $K$. To a $C$-order $A$ in $\Sigma$, M.~Van~den~Bergh assigned in \cite{Michel} a \emph{Brauer-Severi scheme} $\texttt{BS}(A)$ as follows. Let $\trep_nA$ be the scheme of trace preserving representations of $A$, i.e. all $\CC$-algebra morphisms from $A$ to $M_n(\CC)$ that are compatible with the trace maps on both algebras and consider the natural $\GL_n$ action on
$$\trep_nA\times\CC^n$$
given by
$$g.(\phi,v) := (g\phi g^{-1}, gv).$$
By $\texttt{brauer}(A)$ we denote the set of all Brauer-stable points:
$$\texttt{brauer}(A) := \{(\phi,v)\mid \phi(A)v = \CC^n\}.$$
This is also the set of points with trivial stabilizer group, meaning that every $\GL_n$-orbit in $\texttt{brauer}(A)$ is closed and we can form the orbit space 
$$\texttt{BS}(A) = \texttt{brauer}(A)/\GL_n.$$
Moreover, we have a map
$$\pi:\texttt{BS}(A)\twoheadrightarrow\triss_n(A)$$
where $\triss_n(A)$ is the space of isomorphism classes of trace-preserving semi-simple representations (that is, $\triss_n(A)$ is the quotient space of $\trep_n(A)$ under the natural $\GL_n(\CC)$-action) for which $\CC[\triss_n(A)] = C$. The space $\texttt{BS}(A)$ is a projective space bundle over $\triss_n(A)$, but in general not much else is known about these Brauer-Severi schemes. 

In this paper we will deal only with the situation where $A$ is a \emph{noncommutative smooth order} in the sense of \cite{Lieven:NCModel}.
This is equivalent to the condition that $\trep_n(A)$ is a smooth scheme. In this case $\texttt{BS}(A)$ is also a smooth scheme but 
the quotient space $\triss_n(A)$ can contain singularities. This also means we can describe the \'etale local structure of these objects using quivers.

The situation considered above is very similar to that of the construction of desingularizations of singular varieties $V\rightarrow X$ where we have a smooth variety $V$ and a Zariski-open subset of $X$ on which the fibers of the desingularization map consist of a single point. Here, we again have a smooth variety, $\texttt{BS}(A)$ and a Zariski-open subset (the Azumaya locus) on which the fibers are still quite nice, that is, they are projective spaces $\PP^{n-1}$, hence the term hyper-desingularization. In this paper, we wish to determine for which types of central singularities the Azumaya locus actually coincides with the flat locus and, in case the two loci do not coincide, what the fibers of the Brauer-Severi fibration look like over the flat, non-Azumaya points. We provide a complete answer to both questions, the second of which can be seen as a generalization of the work of M.~Artin on the Brauer-Severi fibration over a ramified point, see \cite{Artin}.

The paper is organized as follows: in Section \ref{prelims} we gather together all necessary background information for the remainder of the paper; in Section \ref{flatlocus}, we describe the flat locus of a Brauer-Severi scheme associated to a smooth order; in Section \ref{singularities} we derive all possible central singularities for which the Azumaya locus does not coincide with the flat locus and in Section \ref{fibers} we give a description of the fibers over points in the flat, non-Azumaya locus near such central singularity. Finally, in Section \ref{toric}, we show how these fibers can be seen as toric varieties and we use this to compute their cohomology. 

We want to stress that in the remainder of the paper, whenever Brauer-Severi schemes are mentioned, we will assume that we work over smooth orders.  

\section{Preliminaries}\label{prelims}
We begin by introducing the notions and results we will need throughout the rest of this paper.

\subsection{Definitions and Notations}\label{Definitions}
\begin{defn}\label{QuiverDefinitions}
\begin{itemize}
\item[]
\item A \emph{quiver} is a fourtuple $Q = (Q_0, Q_1, h, t)$ consisting of a set of \emph{vertices} $Q_0$, a set of \emph{arrows} $Q_1$ and two maps $t:Q_1\rightarrow Q_0$ and $h:Q_1\rightarrow Q_0$ assigning to each arrow its tail resp.~its head:
$$\xymatrix@R-1.5pc{\vtx{~} & \vtx{~}\ar[l]_a \\ h(a) & t(a)}.$$
\item A \emph{dimension vector} of a quiver $Q$ is a map $\alpha: Q_0 \rightarrow\NN$ and a \emph{quiver setting} is a couple $(Q,\alpha)$ of a quiver and an associated dimension vector. The dimension vector which is equal to $1$ on all vertices is denoted by $\mathbf{1}$.
\item Fix an ordering of the vertices of $Q$. The \emph{Euler form of a quiver $Q$} is the bilinear form
$$\chi_Q : \NN^{\# Q_0}\times \NN^{\# Q_0} \rightarrow \ZZ$$
defined by the matrix having $\delta_{ij} - \#\{a\in Q_1 \mid h(a) = j, t(a) = i\}$ as element at location $(i,j)$.
\item A quiver is called \emph{strongly connected} if and only if each pair of vertices in its vertex set belongs to an oriented cycle.
\end{itemize}
\end{defn}
A quiver setting is graphically depicted by drawing the quiver and listing in each vertex $v$ either  the dimension $\alpha(v)$, in which case the vertex is encircled, or the name of the vertex itself.

\begin{defn}\label{QuiverRepresentationDefinition}
\begin{itemize}
\item[]
\item An $\alpha$-dimensional representation $V$ of a quiver $Q$ assigns to each vertex $v\in Q_0$ a linear space $\CC^{\alpha(v)}$ and to each arrow $a \in Q_1$ a matrix $V(a) \in M_{\alpha(h(a))\times \alpha(t(a))}(\CC)$. We denote by $\rep(Q,\alpha)$ the space of all $\alpha$-dimensional representations of $Q$. That is,
$$\rep(Q,\alpha) = \bigoplus_{a\in Q_1} M_{\alpha(h(a))\times\alpha(t(a))}(\CC).$$
\item We have a natural action of the reductive group
$$\GL_\alpha := \prod_{v\in Q_0} \GL_{\alpha(v)}(\CC)$$
on a representation $V$ defined by basechange in the vectorspaces. That is
$$(g_v)_{v\in Q_0}.(V(a))_{a\in Q_1} = (g_{h(a)}V(a)g_{t(a)}^{-1})_{a\in Q_1}.$$
\item The quotient space with respect to this action classifies all isomorphism classes of semisimple representations and is denoted by $\iss(Q,\alpha)$. The quotient map with respect to this action will be denoted by 
$$\pi_Q:\rep(Q,\alpha)\twoheadrightarrow\iss(Q,\alpha).$$ 
\item The fibre of $\pi_Q$ in $\pi_Q(0)$ is called the \emph{nullcone} of the quiver setting and is denoted by $\Null(Q,\alpha)$.
\end{itemize}
\end{defn}

\subsection{The \'etale Local Structure of the Brauer-Severi Fibration}
Given a Brauer-Severi scheme $\BS(A)$ over a normal projective variety $X$, the \'etale local structure of the fibration $\BS(A)\twoheadrightarrow X$ near a point $p$ was described in \cite{Lieven:SP}.  The structure of the fiber is completely determined by a triple $(Q_p,\alpha_p,\gamma_p)$ where $(Q_p,\alpha_p)$ is a quiver setting and $\gamma_p\in\NN^{\#(Q_p)_0}$ such that 
$$\sum_{v\in (Q_p)_0}\gamma_p(v)\alpha_p(v) = n.$$
The identification goes as follows (for full details we refer to the papers \cite{Lieven:SP} and \cite{Lieven:NCModel}). Construct a new quiver setting $(\tilde{Q}_p,\tilde{\alpha}_p)$ by adding a vertex $v_0$ to $Q$, adding $\gamma_p(v)$ arrows from $v_0$ to $v$ for each $v\in (Q_p)_0$ and assigning to $v_0$ dimension $1$ and to all other vertices $v$ their original dimension $\alpha(v)$. Now consider the character $\theta_p$ on $\tilde{Q}_p$ which assigns $-n$ to vertex $v_0$ and $\gamma_p(v)$ to vertex $v$ for all $v\in (Q_p)_0$. Consider for this character its space of semistable representations, $\ress_{\theta_p}(\tilde{Q}_p,\tilde\alpha_p)$ in the sense of King \cite{King} then
\begin{theorem}\label{localbrauerstructure}
Given a Brauer-Severi scheme $\BS(A)$ over a normal projective variety $X$ and given a point $p\in X$, we have
$$\pi^{-1}(p) = (\Null(\tilde{Q}_p,\tilde{\alpha}_p)\cap\ress_{\theta_p}(\tilde{Q}_p,\tilde\alpha_p))/\GL_{\tilde{\alpha}_p}$$
Moreover, the dimension vector $\alpha_p$ is such that there exist simple representations in $\rep(Q_p,\alpha_p)$.
\end{theorem}
We will call $(Q_p,\alpha_p)$ the \emph{local quiver setting of $p$} and $(Q_p,\alpha_p,\gamma_p)$ the \emph{local quiver data of $p$}.

\textbf{Remark.} In \cite{Lieven:SP} the results are stated using marked quivers i.e. some of the loops (the marked ones) are only represented 
by traceless matrices. As we are only considering representations in the nullcone (all traces are zero), the distinction between marked and unmarked loops is superfluous.

This theorem reduces the study of the fibers of the Brauer-Severi fibration to the study of moduli spaces of nullcones of quiver settings that have simple representations. Now in \cite{LievenProcesi} a criterium for the existence of simple representations of dimension vector $\alpha$ was given.
\begin{theorem}\label{criteriumforsimple}
Let $(Q,\alpha)$ be a quiver setting such that for all vertices $v$ we have that $\alpha(v)\geq 1$. There exist simple representations of dimension vector $\alpha$ if and only if
\begin{itemize}
\item $Q$ has at least two vertices and is of the form
$$
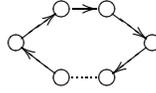
\xymatrix@R=.75pc@C=.75pc{&\vtx{~}\ar[r] & \vtx{~}\ar[dr] \\ \vtx{~}\ar[ur]&  & & \vtx{~}\ar[dl] \\ & \vtx{~}\ar[ul] & \vtx{~}\ar@{.}[l]}$$
with $\alpha = \mathbf{1}$;
\item $Q$ has exactly one vertex, one loop and $\alpha = 1$;
\item none of the above, but $Q$ is strongly connected and
$$\forall v \in Q_0: \chi_Q(\alpha,\varepsilon_v)\leq 0 ~\textrm{and}~\chi_Q(\varepsilon_v,\alpha) \leq 0.$$
Here $\varepsilon_v(w) := \delta_{vw}$ for all $w \in Q_0$.
\end{itemize}
If $\alpha(v) = 0$ for some vertices $v$, $(Q,\alpha)$ has simple representations if $(Q',\alpha')$ has simple representations, where $(Q',\alpha')$ is the quiver obtained by removing all vertices $v$ with $\alpha(v) = 0$.
\end{theorem}

\subsection{Cofree Quiver Settings}
We need one last result before being able to tackle the two questions asked in the introduction: the classification of cofree quiver representations.
\begin{defn}
A quiver setting $(Q,\alpha)$ is called \emph{cofree} if its quotient space $\iss(Q,\alpha)$ is smooth and its nullcone $\Null(Q,\alpha)$ is equidimensional.

A path in a quiver setting will be called \emph{quasiprimitive} if it does not run $n+1$ times through a vertex $v$ with $\alpha(v) = n$. A quasiprimitive path from vertex $v$ to vertex $w$ is depicted as $\xymatrix{v\ar@{~>}[r]&w}$.

A quiver $Q$ is called a \emph{connected sum} of two subquivers $R$ and $S$ in vertex $v$ if $Q_0 = R_0\cup S_0$, $R_0\cap S_0 = \{v\}$, $Q_1 = R_1 \cup S_1$ and $R_1\cap S_1 = \emptyset$.

A quiver setting is called \emph{prime} if it is not a connected sum of two quiver settings in a vertex with dimension $1$, and the prime components of a quiver setting are its maximal prime subquiver settings.
\end{defn}
By \emph{reduction step} $\mathcal{R}^c_I$ we mean the construction of a new quiver from a given quiver by removing a vertex (and connecting all arrows) in the situation illustrated below, where $k$ is not smaller than the number of quasiprimitive cycles through $
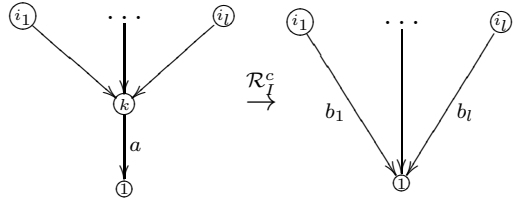
\xymatrix{\vtx{k}}$.
\[\vcenter{\xymatrix{\vtx{i_1}\ar[dr]&\dots\ar[d]&\vtx{i_l}\ar[dl]\\&\vtx{k}\ar[d]^a&\\&\vtx{1}&}}\stackrel{\cR_I^c}{\rightarrow}\vcenter{\xymatrix{\vtx{i_1}\ar[ddr]_{b_1}&\dots\ar[dd]&\vtx{i_l}\ar[ddl]^{b_l}\\&&&\\&\vtx{1}&}}\]
We now have
\begin{theorem}\label{cofree}
A quiver setting $(Q,\alpha)$ is cofree if and only if it can be reduced using $\cR_I^c$ to 
a setting whose prime components are in the list below.
\begin{itemize}
\item[(i)] strongly connected quiver settings $(P,\rho)$ for which
\begin{enumerate}
\item There is a vertex $v \in P_0$  such that $\rho(v) = 1$ and through which all cycles run, 
\item $\forall w\neq v \in P_0: \rho(w) \geq \#\{\xymatrix@=.2cm{v\ar@{~>}[rr]&&w}\}+\#\{\xymatrix@=.2cm{v&& w\ar@{~>}[ll]}\}-1$,
\end{enumerate}
\item[(ii)] quiver settings $(P,\rho)$ of the form
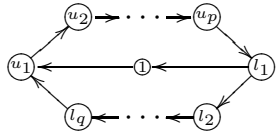
\[
\xymatrix@R=.75pc@C=.75pc{&\vtx{u_2}\ar[r]&\cdots\ar[r]&\vtx{u_p}\ar[dr]&\\\vtx{u_1} \ar[ru]&&\vtx{1}\ar[ll]&& \vtx{l_1} \ar[dl]\ar[ll]\\&\vtx{l_q}\ar[ul]&\cdots \ar[l]&\vtx{l_2}\ar@{->}[l]&}
\]
with $p,q\ge 1$,
such that there is at most one vertex $x$ in the path $\xymatrix{\vtx{u_1} \ar@{~>}[r] & \vtx{l_1}}$ which 
attains the minimal dimension $\min \{u_1,\dots, u_p,l_1,\dots, l_q\}$. 

The special case where $p=0$ and $q\ge 1$ (or vice versa) 
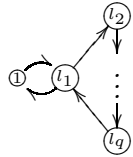
\[
\xymatrix@R=.75pc@C=.75pc{
&&\vtx{l_2}\ar[d]\\
\vtx{1}\ar@/^/[r]&\vtx{l_1}\ar@/^/[l]\ar[ru]&\vdots\ar[d]\\
&&\vtx{l_q}\ar[lu]}
\]
is always cofree.
\item[(iii)] quiver settings of extended Dynkin type $\tilde{A}_n$ with cyclic orientation
\item[(iv)] quiver settings $(P,\rho)$ consisting of two cyclic quivers, $\tilde{A}_{p+s-1}$ and $\tilde{A}_{q+s-1}$, coinciding on $s$ subsequent vertices ($p,q$ can be zero)
$$
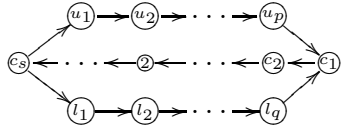
\xymatrix@R=.75pc@C=.75pc{
                                     & \vtx{u_1}\ar[r] & \vtx{u_2}\ar[r] & \dots\ar[r]   &  \vtx{u_p}\ar[dr] & \\
\vtx{c_s}\ar[ur]\ar[dr] & \dots\ar[l]          & \vtx{2}\ar[l]     & \dots\ar[l] & \vtx{c_2}\ar[l] & \vtx{c_1}\ar[l]\\
                                     & \vtx{l_1}\ar[r] & \vtx{l_2}\ar[r] & \dots\ar[r]   &  \vtx{l_q}\ar[ur]
}$$
with
$u_i,l_j\geq 2$ for all $1\leq i \leq p$, $1\leq j \leq q$ and 
all $c_k \ge 4$ except for a unique vertex with dimension $2$.
\end{itemize}
\end{theorem}

\section{The Flat Locus of the Brauer-Severi Fibration}\label{flatlocus}
In this section, we will determine the flat locus of the Brauer-Severi fibration
$$\BS(A)\twoheadrightarrow X,$$
using the quiver description recalled in the Section \ref{prelims}. We have
\begin{theorem}
Let $(Q,\alpha,\gamma)$ be the local quiver data of a point $\xi\in X$, then
$$\dim\pi^{-1}(\xi) = \dim\Null(Q,\alpha) + n - \dim\GL(\alpha).$$
\end{theorem}
\prf Let $\gamma  = (d_1, \dots, d_k)$ and let $\theta$ be the corresponding character for $\tilde{Q}$, then by Theorem \ref{localbrauerstructure} we know that
$$\pi^{-1}(\xi) = (\Null(\tilde{Q},\tilde{\alpha})\cap\rep_\theta(\tilde{Q},\tilde{\alpha}))/\GL(\tilde{\alpha}).$$
Also
$$\Null(\tilde{Q},\tilde{\alpha}) = \Null(Q,\alpha)\times\Affine^{n}$$
because any choice of a representation for an arrow with tail $v_0$ belongs to the nullcone by a straightforward application of the Hilbert criterium. Now for any irreducible component $C$ of $ \Null(Q,\alpha)$ we have that $\rep_\theta(\tilde{Q},\tilde{\alpha})\cap (C\times\Affine^{n}) \neq \emptyset$ only contains $\theta$-stable representations which have stabilizer $\CC^*$ so when $C$ is an irreducible component of maximal dimension we obtain
\begin{eqnarray*}
\dim(\Null(\tilde{Q},\tilde{\alpha})\cap\rep_\theta(\tilde{Q},\tilde{\alpha}))/\GL(\tilde{\alpha})& =& \dim C  + n\\
& &  - (\dim \GL(\tilde{\alpha}) - 1)\\
& =& \dim\Null(Q,\alpha) + n - \dim\GL(\alpha).
\end{eqnarray*}
\eop
This yields
\begin{cor}
A point $\xi\in X$ belongs to the flat locus of the Brauer-Severi fibration of $X$ if and only if its corresponding local quiver setting is cofree.
\end{cor}
\prf
The fibers of the Brauer-Severi fibration over an Azumaya point are isomorphic to $\PP^{n-1}$ so have dimension $n-1$. This means that in order to have minimal dimension we must have $\dim\pi^{-1}(\xi) = n-1$. This is exactly the case when 
$$\dim\Null(Q,\alpha) = \dim\GL(\alpha)-1 = \dim\rep(Q,\alpha)-\dim\iss(Q,\alpha).$$
In combination with the fact that the flat locus of the Brauer-Severi fibration is contained within the smooth locus of the Brauer-Severi fibration by an application of the Popov conjecture for quiver representations (see \cite{Geert:Popov}) we obtain the claim.
\eop
In combination with Theorem \ref{localbrauerstructure}, this means that in order to find the flat locus of the Brauer-Severi fibration of $X$, we have to classify all cofree quiver representations that have simple representations. 
\begin{theorem}\label{th:cofreesimple}
A point $\xi\in X$ belongs to the flat locus of the Brauer-Severi fibration of $X$ if and only if the prime components of its local quiver are
\begin{itemize}
\item either cyclic quiver settings of extended Dynkin type $\tilde{A}$ and dimension vector $\mathbf{1}$;
\item a quiver setting of the form
\[
\xymatrix@R=.75pc@C=.75pc{&\vtx{d}\ar[r]&\cdots\ar[r]&\vtx{d}\ar[dr]&\\
\vtx{d} \ar[ru]&&\vtx{1}\ar[ll]&& \vtx{d} \ar[dl]\ar[ll]\\&\vtx{d-1}\ar[ul]&\cdots \ar[l]&\vtx{d-1}\ar@{->}[l]&}
;
\]
\item a quiver setting of the form
\[
\xymatrix@R=.75pc@C=.75pc{
&&\vtx{d}\ar[d]\\
\vtx{1}\ar@/^/[r]&\vtx{d}\ar@/^/[l]\ar[ru]&\vdots\ar[d]\\
&&\vtx{d}\ar[lu]};
\]
\item a quiver setting of the form
\[
\xymatrix@R=.75pc@C=.75pc{
&&\vtx{d-1}\ar[d]\\
\vtx{1}\ar@/^/[r]&\vtx{d}\ar@/^/[l]\ar[ru]&\vdots\ar[d]\\
&&\vtx{d-1}\ar[lu]}
\]
\item or a connected sum of two cyclic quivers of extended Dynkin type $\tilde{A}$ with dimension vector $2\cdot\mathbf{1}$.
\end{itemize}
\end{theorem}
\prf As $\xi$ belongs to the flat locus of the Brauer-Severi fibration, we know its local quiver setting must have simple representations. Now note that for a dimension vector to have simple representations, the vertex $\xymatrix{\vtx{k}}$ in reduction step $\cR_I^c$ must have dimension $k = 1$. This means there runs only one quasiprimitive path through this vertex, which is only the case if the quiver setting to which the vertex belongs to, is a connected sum of cyclic quivers with dimension vector $\mathbf{1}$. All other cofree quiver settings with simple representations cannot be reduced by reduction step $\cR_I^c$ and hence must have prime components as described in Theorem \ref{cofree}.

The condition $\chi_Q(\alpha,\varepsilon_v)\leq 0$ and $\chi_Q(\varepsilon_v,\alpha)\leq 0$ means that for any given vertex $v$, 
$$\sum_{a\in Q_1, h(a)=v} \alpha(t(a)) \leq \alpha(v) \geq  \sum_{a\in Q_1, t(a)=v} \alpha(h(a)).\quad (*)$$
Consider the prime components of Theorem \ref{cofree} described in $(iv)$. Only a unique $\xymatrix{\vtx{c_i}}$ can have dimension $2$, so $s=1$ by condition $(*)$ and by the same condition $u_i = 2$ and $l_j = 2$ for all $i,j$. This yields a connected sum of cyclic quivers with dimension vector $2.\mathbf{1}$. The component described in $(iii)$ (the cyclic quiver $\tilde{A}$) must obviously have dimension vector $\mathbf{1}$ by Theorem \ref{criteriumforsimple}.  The same argument used for the components listed in (iv) applies to the components described in $(ii)$, which means they only have simple representations when they are of the form
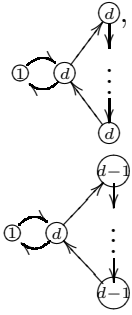
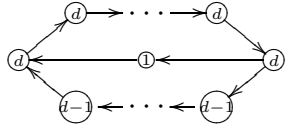
\[
\xymatrix@R=.75pc@C=.75pc{
&&\vtx{d}\ar[d]\\
\vtx{1}\ar@/^/[r]&\vtx{d}\ar@/^/[l]\ar[ru]&\vdots\ar[d]\\
&&\vtx{d}\ar[lu]},
\]
\[
\xymatrix@R=.75pc@C=.75pc{
&&\vtx{d-1}\ar[d]\\
\vtx{1}\ar@/^/[r]&\vtx{d}\ar@/^/[l]\ar[ru]&\vdots\ar[d]\\
&&\vtx{d-1}\ar[lu]}
\]
or
\[
\xymatrix@R=.75pc@C=.75pc{&\vtx{d}\ar[r]&\cdots\ar[r]&\vtx{d}\ar[dr]&\\
\vtx{d} \ar[ru]&&\vtx{1}\ar[ll]&& \vtx{d} \ar[dl]\ar[ll]\\&\vtx{d-1}\ar[ul]&\cdots \ar[l]&\vtx{d-1}\ar@{->}[l]&}
\]

Finally, consider the components described in $(i)$. Let $v$ be a vertex with $\alpha(v) = 1$ through which all cycles run. Because of this condition, there is a vertex $w_0$ that has only incoming arrows from $v$. Should $\alpha(w_0) > 1$, the existence of simple representations implies there are $\alpha(w_0)$ arrows from $v$ to $w$. The second condition on $w$ then implies there is only one arrow leaving, and this to a vertex $w_1$ with dimension $\alpha(w_1) = \alpha(w)$. As all cycles run through $v$, this vertex also has only one arrow leaving to a vertex $w_2$ with dimension $\alpha(w_2) = \alpha(w)$. This argument can be repeated until we find a vertex $w_k$ that is connected to $v$ with exactly one arrow, so all vertices $w_i$ for $0\leq i \leq k$ must have $\alpha(w_i) = 1$. This means the quiver setting has a prime component that is a cyclic quiver setting with dimension vector $\mathbf{1}$. Removing this prime component then yields by induction on the number of vertices that the quiver setting we started with was a connected sum of cyclic quivers with dimension vector $\mathbf{1}$.
\eop

\section{Central Singularities and the Flat Locus}\label{singularities}
In this section, we give a classification of all reduced (in the sense of \cite{SOS}) central singularities where the flat locus of the Brauer-Severi fibration does not coincide with the Azumaya locus. This translates in finding all reduced non cofree quiver settings $(Q,\alpha)$, which have a local quiver setting that is cofree. Recall from \cite{LievenProcesi} the construction of a local quiver setting $(Q_p,\alpha_p)$ from a quiver $(Q,\alpha)$:
for a given semisimple representation 
\[
S_1^{\oplus a_1}\oplus\ldots\oplus S_k^{\oplus a_k}
\]
of the quiver $(Q,\alpha)$, where $S_i$ has dimension vector $a_i$ (and a combinatorial property of simplicity described in Theorem \ref{criteriumforsimple}), we construct the local quiver $(Q_p,\alpha_p)$ by the following data:
\begin{itemize}
\item $(Q_p)_0$ consist of $k$ vertices.
\item The number of arrows between verices $u$ and $v$ is given by $\delta_{uv}-\chi_Q(\alpha_u,\alpha_v)$.
\item $\alpha_p:=(a_1,\ldots,a_k)$.
\end{itemize} 

In the remainder of this section we will use this combinatorial description to obtain the desired classification. From now on, we will call a quiver setting simple if it has simple representations. 

A characterization of reduced quiver settings (with at least 2 vertices) is the following:
\begin{defn}\label{def:reduced}
A strongly connected quiver setting $(Q,\alpha)$ with at least 2 vertices is called reduced if and only if
\begin{itemize}
\item every vertex $t$ with no loop has $\chi_Q(\alpha,\varepsilon_t)\leq -1$ and $\chi_Q(\varepsilon_t,\alpha)\leq -1$.
\item every vertex $t$ with one loop has $\chi_Q(\alpha,\varepsilon_t)\leq -2$ and $\chi_Q(\varepsilon_t,\alpha)\leq -2$.
\item all vertices with dimension 1 do not have loops.
\end{itemize}
\end{defn}

\begin{lemma}\label{lem:twovertices}
A non cofree quiver setting with a cofree local quiver has at least 2 vertices.
\end{lemma}
\prf
From \cite{Raf_Geert:Cofree} we easily deduce that a non cofree quiver setting with one vertex has at least 2 loops and dimension at least 3 or at least 3 loops and dimension at least 2. In both cases one easily verifies that too many arrows appear in the local quivers.
\eop
From now on we suppose $(Q,\alpha)$ to be a reduced quiver setting with at least 2 vertices (this implies that $(Q,\alpha)$ is not cofree). We first determine the possible local quiver settings.

\begin{lemma}\label{lem:twoloops}
A cofree local quiver setting of $(Q,\alpha)$ can never have 2 or more loops at a vertex $v$ with dimension $\alpha(v) = k \geq 2$.
\end{lemma}
\prf
The only situation of Theorem \ref{th:cofreesimple} where two loops occur is
\begin{equation}\label{eq:quiver1}
\xymatrix{\vtx{2} \ar@(ld,lu) \ar@(rd,ru)}
\end{equation}
Suppose one can find a reduced quiver setting $(Q,\alpha)$ with $n$ vertices such that it has (\ref{eq:quiver1}) as a local quiver. Then we know that 2 divides $\alpha_i$ for all $i$. Since $(Q,\alpha)$ is reduced, $-\chi(\alpha,\varepsilon_i)\geq 1$ and of course also $-\chi(\alpha/2,\varepsilon_i)\geq 1$. But then the number of loops of the local quiver is
\[
1-\chi_Q(\frac{\alpha}{2},\frac{\alpha}{2})=1-\sum_{i=1}^n\frac{\alpha}{2}\chi_Q(\frac{\alpha}{2},\varepsilon_i) \geq 3.
\]
\eop

\begin{lemma}\label{lem:an}
$(Q,\alpha)$ is simple and $\chi_{Q}(\alpha,\alpha)=0$ if and only if $Q=\tilde{A}_n$ and $\alpha=(1,\ldots,1)$.
\end{lemma}
\prf
By linearity of $\chi_Q$, it is clear that $\chi_{Q}(\alpha,\alpha)=0$ is equivalent to $\chi_{Q}(\alpha,\varepsilon_t)=0$ for all $t$. Looking at the vertex with minimal dimension, one easily deduce that all dimensions must be the same. The only strongly connected quiver allowing this is an $\tilde{A}_n$. The simplicity condition gives us dimension vector $(1,\ldots,1)$.
\eop

\begin{theorem}
Suppose $(Q,\alpha)$ has a cofree local quiver setting $(Q_p,\alpha_p)$, then $(Q_p,\alpha_p)$ is always a connected sum of $\tilde{A}_i$, where every vertex has dimension $1$ (and a number of loops).
\end{theorem}
\prf Suppose we have a cofree local quiver with a vertex of dimension $k \geq 2$. This vertex $v$ corresponds to a simple quiver setting $(Q_v,\alpha_v)$, with $Q_v$ a subquiver of $Q$ and $\alpha_v$ a sub-dimension vector of $\alpha$, nonzero on $Q_v$ and zero elsewhere. From Lemma \ref{lem:twoloops} and the fact that the number of loops in vertex $v$, is given by
\[
1-\chi_Q(\alpha_u,\alpha_u)
\]
and therefore is at least 1, we know that vertex $v$ has exactly one loop. This also implies by Lemma \ref{lem:an} that $(Q_u,\alpha_u)=(\tilde{A}_n,\bf{1})$. The number of non-loop arrows arriving in vertex $v$ is given by
\begin{align*}
&\sum_{w \neq v} -\chi_Q(\alpha_w,\alpha_v) \\
&= - \chi_Q(\alpha-k\alpha_v,\alpha_v)\\
&= - \chi_Q(\alpha,\alpha_v)+k\underbrace{\chi_Q(\alpha_v,\alpha_v)}_{=0}\\
&= - \sum_{t \in (Q_v)_0}\chi_Q(\alpha,\varepsilon_t)\alpha_v(t)\geq 2
\end{align*}
where the last inequality holds because $(Q,\alpha)$ is reduced. A situation with one loop and 2 more incoming arrows in a vertex with dimension greater than $k$ can (Theorem \ref{th:cofreesimple}) never be cofree.

We find that all vertices of the local quiver setting need to have dimension 1. The only possibilities left are, according Theorem \ref{th:cofreesimple}, connected sums of $\tilde{A}_i$ having dimension 1 on each vertex.
\eop

Now that we have found all possible local quiver settings, let us look at all possible reduced quiver settings that have a connected sum of $\tilde{A}_i$, with dimension vector $\mathbf{1}$ as local quiver.

\begin{theorem}
The reduced quiver settings that have a connected sum of $\tilde{A}_i$, with dimension vector $\mathbf{1}$ as local quiver, look the same as their local quiver, with every vertex replaced by a simple subquiver, where the $\tilde{A}_i$-arrows always end and start in a vertex with dimension 1.
\end{theorem}
We illustrate this theorem with an example:
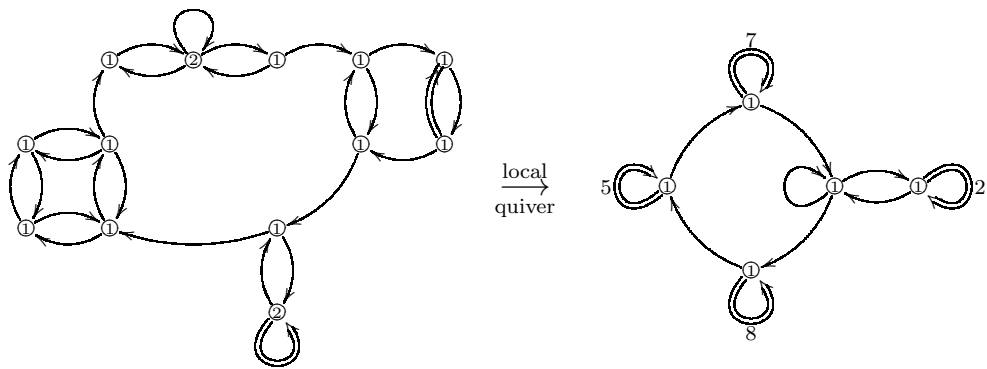
\[
\vcenter{\xymatrix{ 
& \vtx{1} \ar@/^/[r]& \vtx{2} \ar@(lu,ru)\ar@/^/[r]\ar@/^/[l] & \vtx{1} \ar@/^/[l] \ar @/^/[r] &   \vtx{1}\ar@/^/[d] \ar @/^/[r] &  \vtx{1}\ar@/^/[d] \\
 \vtx{1}\ar@/^/[r]  \ar@/^/[d]& \vtx{1}\ar @/^/[u]\ar@/^/[d] \ar@/^/[l] & & & \vtx{1}\ar@/^/[u]\ar @/^/[dl] & \vtx{1}\ar@/^/@{=>}[u] \ar @/^/[l] \\
\vtx{1}\ar@/^/[u]  \ar@/^/[r] & \vtx{1}\ar@/^/[u]  \ar@/^/[l] & & \vtx{1}\ar@/^/[d] \ar @/^/[ll]\\
& & & \vtx{2}\ar@/^/[u]\ar@{=>}@(ld,rd) 
}
}\quad \underset{\text{quiver}}{\overset{\text{local}}{\longrightarrow}} \quad 
\vcenter{
\xymatrix{
& \vtx{1} \ar@/^/[dr] \ar@(ul,ur)@{=>}^7\\
\vtx{1} \ar@(dl,ul)@{=>}^5 \ar@/^/[ur] && \vtx{1} \ar@(dl,ul) \ar@/^/[r]\ar@/^/[dl] & \vtx{1} \ar@/^/[l] \ar@(ur,dr)@{=>}^2\\
& \vtx{1} \ar@/^/[lu] \ar@(dl,dr)@{=>}_8
}}
\]

\prf
We first show that the quivers described in the theorem have a local quiver that is a connected sum of $\tilde{A}_i$ with dimension 1 on every vertex.

Assume $(Q,\alpha)$ of the form 
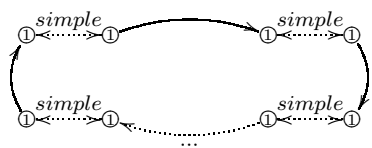
\[ 
\xymatrix{
\vtx{1} \ar@{<..>}^{simple}[r]& \vtx{1} \ar@/^/[rr] & & \vtx{1} \ar@{<..>}^{simple}[r]& \vtx{1} \ar@/^/[d]\\
\vtx{1} \ar@/^/[u] \ar@{<..>}^{simple}[r] & \vtx{1} & & \vtx{1} \ar@/^/@{..>}^{\cdots}[ll] \ar@{<..>}^{simple}[r]& \vtx{1}
}
\]
with $\alpha=(\alpha_1, \alpha_2, \ldots, \alpha_n)$ where $\alpha_j$ is the dimension vector of the simple subquiver placed at vertex $j$ of $\tilde{A}_n$.

The Euler form of this quiver setting is given by the blockmatrix
$$
\begin{pmatrix}
\fbox{$A_1$} & \fbox{$\epsilon_{21}$} & 0 & 0 & \cdots & 0\\
0 & \fbox{$A_2$} & \fbox{$\epsilon_{32}$} & 0 & \cdots & 0\\
0 & 0 & \fbox{$A_3$} & \fbox{$\epsilon_{43}$} & \cdots & 0\\
0 & 0 & 0 & \fbox{$A_4$} & \ddots & 0\\
\vdots & \vdots & \vdots  & & \ddots & \fbox{$\epsilon_{n,n-1}$}\\
\fbox{$\epsilon_{1n}$} & 0 & 0 & 0 & 0 & \fbox{$A_n$}
\end{pmatrix}
$$
where $A_j$ consist of the Euler matrix of the simple quiver at vertex $j$ of $\tilde{A}_n$ and $\epsilon_{ij}$ is a zero matrix with a -1 at exactly one entry.

If we construct the local quiver by splitting the dimension vector $\alpha$ in $n$ components
$$
(\alpha_1,0,\ldots,0)\oplus (0,\alpha_2,0,\ldots,0)\oplus \ldots \oplus (0,\ldots,0,\alpha_n).
$$
we get a local quiver that has $n$ vertices and with exception of loops, the only arrows occuring are exact the same as in the $\tilde{A}_n$. This results in 
$$
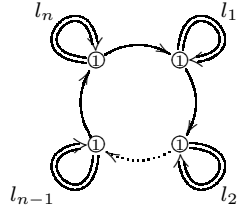
\xymatrix{
\vtx{1} \ar@/^/[r] \ar@{=>}@(l,u)^{l_n} & \vtx{1} \ar@/^/[d] \ar@{=>}@(u,r)^{l_1}\\
\vtx{1} \ar@/^/[u] \ar@{=>}@(d,l)^{l_{n-1}} & \vtx{1} \ar@/^/@{..>}[l] \ar@{=>}@(r,d)^{l_2}
}
$$
which is case 1 of the cofree simple quivers of Theorem \ref{th:cofreesimple}.

More generally, we start with a connected sum of $\tilde{A}_i$, where every vertex is replaced by a simple subquiver and  $\tilde{A_i}$-arrows always start and end in a vertex with dimension 1. We construct a local quiver in the same way as above. The simple subquivers will be packed in 1 vertex and for the same reason as above, we get a local quiver that looks exactly the same as the original,  replacing the simple subquivers by vertices with dimension 1, with a number of loops. We find a connected sum of cyclic quivers.

Remains to show that these quiver settings are the only possible reduced quiver settings that have a cofree local quiver. This result is an immediate consequence of the fact that it is impossible for 2 simple components of the local quiver to have a common vertex with nonzero dimension.

To see that it is impossible to have a common vertex, we use, as before, only reduced quiver settings $(Q,\alpha)$ with at least 2 vertices. Suppose $(Q,\alpha)$ has a local quiver with at least 2 components that overlap in a vertex $t \in Q_0$ (there exist components $u$ and $v$ with $\alpha_u(t),\alpha_v(t)\geq 1$). We now shall show that the local quiver $(Q_p,\alpha_p)$ always has a subquiver of the form\\
\parbox{1cm}{(i)} \parbox{3cm}{$
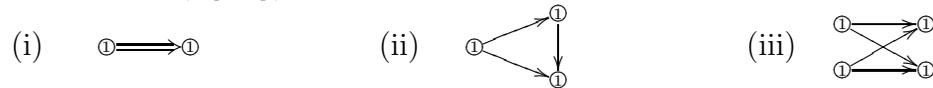
\xymatrix{ \vtx{1} \ar@{=>}[r] & \vtx{1} }$}
\hfill
\parbox{1cm}{(ii)} \parbox{3cm}{ $ \xymatrix@R=.40pc{ & \vtx{1}\ar[dd] \\ \vtx{1} \ar[ru] \ar[rd] \\ & \vtx{1} }$}
\hfill
\parbox{1cm}{(iii)} \parbox{3cm}{$ \xymatrix@R=.80pc{ \vtx{1} \ar[r] \ar[rd] & \vtx{1} \\ \vtx{1} \ar[r] \ar[ru] & \vtx{1}} $} \\
(where in (iii) the two left vertices may coincide) and therefore the local quiver cannot be a connected sum of $\tilde{A}_i$'s.

Let us now define for a set of vertices of the local quiver, say $\Omega$, the function 
\[
\kappa_t(\Omega):=\sum_{u \in \Omega} -\chi_Q(\alpha_u,\varepsilon_t).
\]
If we take $\Omega=(Q_p)_0$, by Definition \ref{def:reduced}, $\kappa_t(\Omega) \geq 1$ and even $\kappa_t(\Omega) \geq 2$ if there is at least 1 loop in $t$. If we now remove vertices $u$ from $\Omega$ with $-\chi_Q(\alpha_u,\varepsilon_t) \leq 0$, $\kappa_t(\Omega)$ still has the same lower bound.

Our next objective is to remove vertices from $\Omega$ in such a way that $\Omega$ has exactly 2 vertices $u_i$ with $\alpha_{u_i}(t)>0$.

We start with $\Omega=(Q_p)_0$ and remove vertices $u$ with $-\chi_Q(\alpha_u,\varepsilon_t) \leq 0$ and $\alpha_u(t)>0$ from $\Omega$, but we stop removing vertices if there are only 2 components $u_i$ with $\alpha_{u_i}(t)>0$ left in $\Omega$. If it is not possible to reach an $\Omega$ with 2 components $u_i$ with $\alpha_{u_i}(t)>0$ by this procedure, we just remove vertices until no vertex $w$ with $\chi_Q(\alpha_w,\varepsilon_t) \leq 0$ and $\alpha_w(t)>0$ exists in $\Omega$. We call the resulting set $\Omega_t$.

If $\Omega_t$ has at least 3 components $u_i$ with $\alpha_{u_i}(t)>0$, each of these $u_i$ has $-\chi_Q(\alpha_{u_i},\varepsilon_t) \geq 1$, and the number of arrows between $u_i$ and $u_j$ is given by
\[
-\chi_Q(\alpha_{u_i},\alpha_{u_j})=\sum_{s \in Q_0} -\chi_Q(\alpha_{u_i},\varepsilon_s)\alpha_{u_j}(s)\geq -\chi_Q(\alpha_{u_i},\varepsilon_t)\alpha_{u_j}(t) \geq 1
\]
where the first inequality is obtained by simplicity of the components of the local quiver. This result in a subquiver of $Q_p$ looking like (ii).

If $\Omega_t$ has 2 such components $u$ and $v$, three situations may occur, but they all will lead to a situation (i), (ii) or (iii):

{\bf A.} $(Q_u)_0 \not\subset (Q_v)_0$ and $(Q_v)_0 \not\subset (Q_u)_0$. The number of arrows from $u$ to $v$ is given by
\begin{equation}\label{eq:numarrows}
-\sum_{s \in (Q_u)_0 \cap (Q_v)_0}\chi_Q(\alpha_u,\varepsilon_s)\alpha_v(s)-\sum_{s \in (Q_v)_0 \backslash (Q_u)_0}\chi_Q(\alpha_u,\varepsilon_s)\alpha_v(s)
\end{equation}
By simplicity of $(Q_u,\alpha_u)$, the first sum is always $\geq 0$, because there exist vertices in $(Q_v)_0 \backslash (Q_u)_0$, the second sum is always at least 1. So we have at least 1 arrow from $u$ to $v$ and by the same argument, we have an arrow in the other direction. Since $\kappa_t(\Omega_t)\geq 1$, we have 
\begin{equation}\label{eq:kappa}
-\chi_Q(\alpha_u,\varepsilon_t)-\chi_Q(\alpha_v,\varepsilon_t)+\sum_{w \in \Omega_t\backslash \{u,v\}}\chi_Q(\alpha_w,\varepsilon_t) \geq  1
\end{equation}
If the first term of (\ref{eq:kappa}) is at least one we get an extra arrow from $u$ to $v$ due to the first sum of (\ref{eq:numarrows}), so we are left with $\xymatrix@C=1.10pc{u \ar@{=>}[r] &v}$ in our local quiver. If the second term is at least one we get $\xymatrix@C=1.10pc{v \ar@{=>}[r] &u}$ in our local quiver. If the third term is at least one, we get the situation (ii). 

{\bf B.} $(Q_u)_0 \subsetneq (Q_v)_0$. If $Q_u$ has exactly one vertex, the number of vertices between $v$ and $u$ is given by $-\chi_Q(\alpha_v,\varepsilon_t)\alpha_u(t)$ and the condition on $\kappa(\Omega_t)$ translates in 
\begin{equation*}
-\chi_Q(\alpha_v,\varepsilon_t)+\sum_{w \in \Omega_t\backslash \{u,v\}}\chi_Q(\alpha_w,\varepsilon_t) \geq  2
\end{equation*}
(this is independent of the number of loops in $t$) and we see that all possible situations lead to a subquiver of $(Q_p)$ of the form (i), (ii) or (iii).

If $Q_u$ has  at least 2 vertices, we can look again at equation (\ref{eq:numarrows}) from situation A, for arrows from $u$ to $v$. We also use (\ref{eq:kappa}), to see we get another arrow from $u$ tot $v$ or vice versa, or an arrow to $u$ and $v$ from another vertex $w$. Let us now look at another vertex $t_1$ of $(Q_u)_0\cap (Q_v)_0$. If here too we have $\kappa_{t_1}(\Omega_t)\geq 1$ this leads to another arrow from $v$ to $u$. If $\kappa_{t_1}(\Omega_t)=0$, there exist one vertex $w'$ of the local quiver, which we have removed while constructing $\Omega_t$. This vertex has $-\chi_Q(\alpha_{w'},\varepsilon_{t_1})$ and this leads to an arrow from $w'$ to both $u$ and $v$. Again we have a subquiver of $(Q_p)$ of the form (i), (ii) or (iii).

{\bf C.} $(Q_u)_0 = (Q_v)_0$. If $(Q_u)_0$ has 1 vertex, the same argument as in B holds. Suppose $(Q_u)_0$ has at least 2 vertices. Again, relation (\ref{eq:kappa}) gives us another arrow from $u$ tot $v$ or reverse, or an arrow to $u$ and $v$ from another vertex $w$. Using the same stategy as in B, we get, for another vertex $t_1$ in the intersection, another arrow. If the intersection consist of 3 or more vertices, we get situation (i) - (iii). The last possibility to consider is the case where we have 2 vertices $t,t_1$ in the intersection. By the same reasoning as before $t$ and $t_1$ each give one arrow $u$ to $v$ or a reverse arrow, or an arrow from $w$ to $u$ and $v$. All cases gives (i) to (iii), except when we have exactly one arrow from $u$ to $v$ and vice versa. A small argument shows that this too is impossible, because this yields $-\chi_Q(\alpha_u,\varepsilon_t)\geq 1$, $-\chi_Q(\alpha_v,\varepsilon_t)=0$, $-\chi_Q(\alpha_v,\varepsilon_{t_1})\geq 1$ and $-\chi_Q(\alpha_v,\varepsilon_{t_1})=0$ and the only way to obtain this is that $(Q_u,\alpha_u)$ or $(Q_v,\alpha_v)$ are not simple quiver settings.
\eop

\section{The Brauer-Severi Fibration over the Flat Locus}\label{fibers}
From the previous section, we know the only flat, non-Azumaya settings that can occur near a central singularity are cyclic quiver settings with dimension vector $\mathbf{1}$. In this section, we give a description of the fibers of the Brauer-Severi fibration over these points in the flat locus.
\begin{lemma}
Let $(Q,\alpha)$ be a quiver setting that is a connected sum of $k$ cyclic quivers $\tilde{A}_{n_i}$, $1\leq i \leq k$ with dimension vector $\alpha = \mathbf{1}$, 
$$Q = \tilde{A}_{n_1}\#_{v_1}\dots\#_{v_{k-1}} \tilde{A}_{n_k} $$
then $\Null(Q,\alpha)$ has $(n_1+1).\dots.(n_k+1)$ irreducible components $T_i$, each of which is a tree that is a connected sum of $k$ quivers of type $A_{n_i+1}$.
\end{lemma}
\prf
Indeed, we know that a representation $V$ lies in the nullcone if and only if all traces along oriented cycles in $Q$ become zero. This condition is equivalent to choosing at least one arrow in each component that gets assigned to $0$ by $V$. But then choosing one arrow in each component gives a closed irreducible subset of the nullcone, and permuting the arrow chosen to be zero yields a covering of the nullcone by irreducible subsets, none of which lies in the union of the others.
\eop
Now let $\xi\in\Flat(Q,\alpha)$ with a local quiver setting $(Q_\xi,\mathbf{1})$ that is a connected sum of cyclic quivers and with $\gamma_\xi = (d_1,\dots, d_X)$. By Theorem \ref{localbrauerstructure} we know that each irreducible component $C$ of $\pi^{-1}(\xi)$ is described by the moduli space $\moss_\theta(\tilde{T},\mathbf{1})$ where $T$ is an irreducible component of of $\Null(Q_\xi,\mathbf{1})$ and hence by the previous lemma a tree. In the remainder of this section we will describe these components. In order to do so, we need some additional definitions.
\begin{defn}
Let $(T,\mathbf{1})$ be a quiver setting of which the underlying quiver is a tree.
\begin{enumerate}
\item A vertex $v \in T_0$ is called a \emph{root vertex} if it is a sink, that is, there are no arrows $a\in T_1$ such that $t(a) = v$.
\item A \emph{rooted tree} is a quiver $T$ that is a tree and for which there exists a unique root vertex.
\item For a root vertex $v$, we denote by $T(v)$ the maximal rooted subtree of $T$ with root vertex $v$.
\item For a vertex $w$ in $T(v)$, we let the \emph{root distance to $v$} be the length of the path connecting $w$ to $v$, and denote this by $\mathcal{D}_v(w)$.
\item For a rooted subtree $T(v)$ we let $h_v := \max\{\mathcal{D}_v(w)\mid w\in T(v)_0\}$ and call this the \emph{height} of the subtree $T(v)$.
\end{enumerate}
\end{defn}
Now let $v$ be a root vertex in the irreducible component $T$ of $\Null(Q_\xi,\mathbf{1})$, then we will assign a graph to $T(v)$ as follows. Assign to each vertex $w\in T(v)_0$ the projective space $\PP^{N_w}$ with 
$$N_w = \sum_{\xymatrix{\vtx{w} & \vtx{u}\ar@{~>}[l]}} d_u -1.$$
For each vertex $w$ in $T(v)$ fix an ordering on the arrows entering $w$, denoting them $a_1^w, \dots a_{r_w}^w$. Denote the vertices at root distance $s$ by $w^s_i$ with $i$ numbering the vertices at the fixed root distance $s$, grouping tails of arrows with the same head together. We now construct a series of projection maps
$$
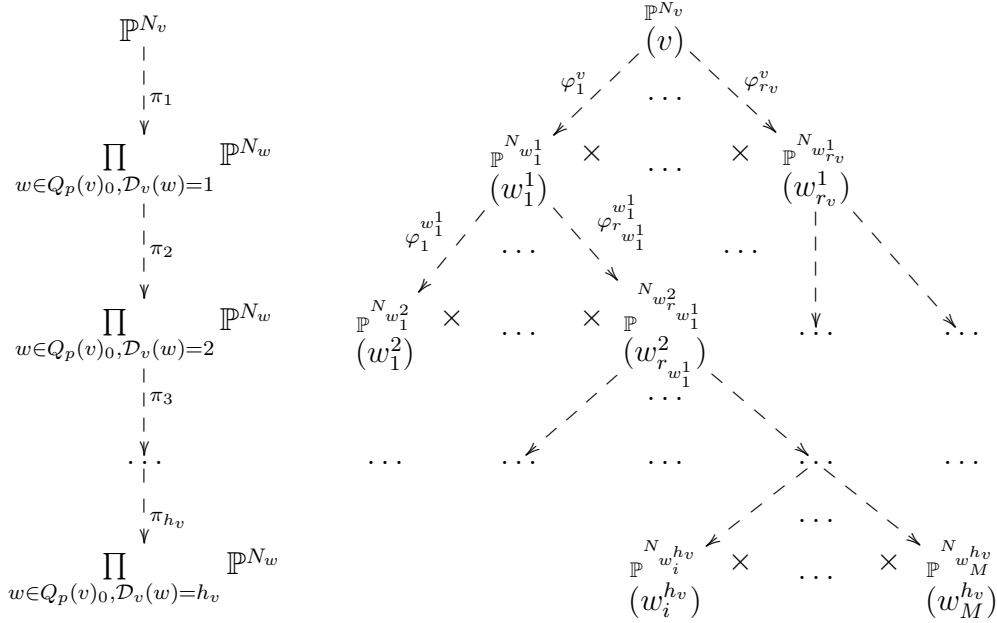
\xymatrix{
\PP^{N_v}\ar@{-->}[d]^{\pi_1} 
& & & \stackrel{\PP^{N_v}}{(v)}
\ar@{-->}[dl]_{\varphi_1^v}\ar@{}[d]|{\displaystyle\dots}\ar@{-->}[dr]^{\varphi_{r_v}^v}
\\
\prod\limits_{w\in Q_p(v)_0, \mathcal{D}_v(w) = 1} \PP^{N_w}\ar@{-->}[d]^{\pi_2} 
& &
\stackrel{\PP^{N_{w_{1}^1}}}{(w_1^1)}
\ar@{}[r]^{\displaystyle\times}
\ar@{}[d]|{\displaystyle\dots}
\ar@{-->}[dl]_{\varphi_1^{w^1_1}} 
\ar@{-->}[dr]^{\varphi_{r_{w_1^1}}^{w^1_1}}
& \dots\ar@{}[r]^{\displaystyle\times}
& 
\stackrel{\PP^{N_{w_{r_v}^1}}}{(w_{r_v}^1)}
\ar@{-->}[d]
\ar@{-->}[dr]
\ar@{}[dl]|{\displaystyle\dots}
\\
\prod\limits_{w\in Q_p(v)_0, \mathcal{D}_v(w) = 2} \PP^{N_w}\ar@{-->}[d] ^{\pi_3}
 &
\stackrel{\PP^{N_{w^2_1}}}{(w^2_1)}
\ar@{}[r]^{\displaystyle\times}
&
\dots 
\ar@{}[r]^{\displaystyle\times}
&
\stackrel{\PP^{N_{w^2_{r_{w_1^1}}}}}{(w^2_{r_{w_1^1}})}
\ar@{-->}[dr]
\ar@{}[d]|{\displaystyle\dots}
\ar@{-->}[dl] 
& \dots 
& \dots
\\
\dots\ar@{-->}[d]^{\pi_{h_v}}  
 & \dots&\dots &\dots & \dots\ar@{-->}[dr]\ar@{}[d]|{\displaystyle\dots}\ar@{-->}[dl] & \dots
\\
\prod\limits_{w\in Q_p(v)_0, \mathcal{D}_v(w) = h_v} \PP^{N_w} 
& & &
\stackrel{\PP^{N_{w^{h_v}_i}}}{(w^{h_v}_i)}
\ar@{}[r]^{\displaystyle\times}
 &
\dots
 \ar@{}[r]^{\displaystyle\times}
 & 
\stackrel{\PP^{N_{w^{h_v}_{M}}}}{(w^{h_v}_M)}
}$$
where 
$$\pi_i = \prod\limits_{w\in Q_p(v)_0, \mathcal{D}_v(w) = i-1} \pi_i^w$$
with
$$\pi_i^w = \prod_{j=1}^{r_w} \varphi^w_j$$
where $\varphi^w_j$ is the projection on the projective coordinates numbered from $1 + \sum_{k=1}^{j-1}(N_{a^w_k}+1)$ to $\sum_{k=1}^{j}(N_{a^w_k}+1)$.
We will denote the graph of this collection of rational maps by $\Graph(v)$.

Let $v$ and $w$ be two root vertices, then their rooted subtrees must have a common subquiver, denoted by $T(v)\cap T(w)$ which is again a rooted tree. Denote the root vertex of this tree by $v\cap w$. We then have
\begin{theorem}\label{FibersAreGraphs}
Let $\xi\in\Flat(Q,\alpha)$ with local quiver data $(Q,\mathbf{1},\gamma)$ where $Q$ is a connected sum of cyclic quivers. Let $(T,\mathbf{1})$ be the quiver setting of an irreducible component of $\Null(Q,\mathbf{1})$. Now let the rooted subtrees of $T$ be connected as follows (we list the roots of the rooted subtrees):
$$
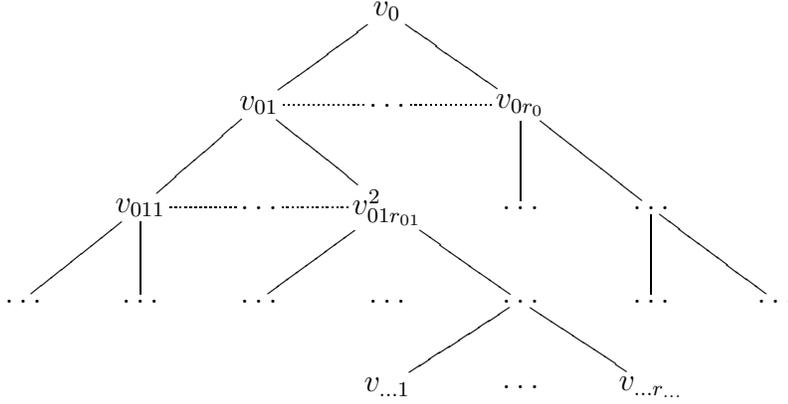
\xymatrix{
& & & v_0\ar@{-}[dl]\ar@{-}[dr] & & & \\
& &v_{01}\ar@{.}[r]\ar@{-}[dl]\ar@{-}[dr] & \dots\ar@{.}[r] & v_{0r_0}\ar@{-}[d]\ar@{-}[dr] & & \\
& v_{011}\ar@{.}[r]\ar@{-}[dl]\ar@{-}[d] & \dots\ar@{.}[r] & v^2_{01r_{01}}\ar@{-}[dl]\ar@{-}[dr] & \dots& \dots\ar@{-}[dr]\ar@{-}[d]& \\
\dots &\dots  &\dots &\dots &\dots\ar@{-}[dl]\ar@{-}[dr] &\dots &\dots  \\
& & & v_{\dots 1} &\dots & v_{\dots r_{\dots}}& 
}$$
That is, $T(v_0)$ has common subgraphs with $T(v_{0i})$ for $1\leq i \leq r_0$ but not with any other $T(w)$ for $w\not\in\{v_{01},\dots, v_{0r_0}\}$, likewise for $T(v_{01})$, and so on. Now let
$$\mathcal{F}_0 = (\dots ((\Graph(v_0)\times_{\Graph(v_0\cap v_{01})}\Graph(v_{01}))\times_{\Graph(v_0\cap v_{02})} \Graph(v_{02})\dots)\times_{\Graph(v_0\cap v_{0r_0})}\Graph(v_{0r_0})),$$
$$\mathcal{F}_1 = (\dots ((\mathcal{F}_0\times_{\Graph(v_{01}\cap v_{011})}\Graph(v_{011}))\times_{\Graph(v_{01}\cap v_{012})} \Graph(v_{012})\dots)\times_{\Graph(v_{0r_0}\cap v_{0r_0r_{0r_0}})}\Graph(v_{0r_0r_{0r_0}})),$$
$$\dots$$
\se{\mathcal{F}_n = (\dots ((\mathcal{F}_{n-1}\times_{\Graph(v_{01\dots 1}\cap v_{01\dots 11})}\Graph(v_{01\dots 11}))\times_{\Graph(v_{01\dots 1}\cap v_{01\dots 12})} \Graph(v_{01\dots 12})\dots)\\
\times_{\Graph(v_{0r_0}\cap v_{0r_0r_{0r_0}\dots r_{0r_0r_{0r_0}\dots}})}\Graph(v_{0r_0r_{0r_0}\dots r_{0r_0r_{0r_0}\dots}})),
}
then the irreducible component of $\pi^{-1}(p)$ corresponding to $T$ is equal to
$\mathcal{F}_n$.
\end{theorem}
\prf
We must describe $\moss_\theta(\tilde{T},\tilde{\mathbf{1}})$. We will first describe $\ress_\theta(\tilde{T},\tilde{\mathbf{1}})$ and then see what the action of $\GL(\tilde{\mathbf{1}})$ on this subspace of semistable representations does. First of all, let $v$ be a root vertex  in $T$  and let $\tilde{T}(v)$ be the subquiver in $\tilde{T}$ corresponding to the rooted subtree of $v$. It is obvious that a representation $V$ is semistable for this setting if and only if each vertex in $T$ is reached by $V$. For a given vertex $w \in T(v)$ this means that the representation must be non-zero along at least one of the paths from $v_0$ to $w$, and this for each $w$. For any vertex $w$, denote the arrows from $v_0$ to $w$ by $x^w_1, \dots, x^w_{\gamma_w}$. By abuse of notation, we will use the same notation for both the arrow $a$ and the value assigned to it by $V$. Now define for a top vertex $t$
$$\mathbf{t} := (x^t_1, \dots, x^t_{\gamma_t}) \in \CC^{\gamma_t}$$
and define inductively for any vertex $w$ with incoming arrows $a_1, \dots, a_{i_w}$
$$\mathbf{w} := (x^w_1, \dots, x^w_{\gamma_w},a_1\mathbf{t(a_1)}, \dots, a_{i_w}\mathbf{t(a_{i_w})}) \in \CC^{N_w+1}.$$
Then the semistability condition for $V$ yields that $\mathbf{w} \neq 0$ for all $w$. Moreover, we may identify $V$ with its image in $\prod_{w\in TT_0}\CC^{N_w+1}$ under the map $V\mapsto (\mathbf{w})_{w\in T_0}$. The action of $\GL(\tilde{\textbf{1}})$ on $V$ translates in the natural action of $\prod_{w\in T_0} \CC^*$ on $\prod_{w\in T_0}\CC^{D_w+1}$ by left multiplication. This means the orbit of $\mathbf{w}$ corresponds to a point in $\prod_{w\in T_0}(\PP)^{N_w}$. Denote by $$\overline{\mathbf{w}} \in \PP^{N_w}$$
the projective coordinates obtained from $\mathbf{w}$, then the orbit of $V$ is a $M$-tuple of projective points (with $M$ the number of vertices in $T$):
$$
\mathcal{O}_V = 
(\dots, \overline{\mathbf{v}}, \overline{\mathbf{t(a^v_1)}},\dots, \overline{\mathbf{t(a^v_{n_v})}}, \dots, \overline{\mathbf{w}}, \overline{\mathbf{t_1}},\dots, \overline{\mathbf{t_{n_w}}},\dots).
$$
For the rooted subtree with root vertex $v$, the points corresponding to the rooted subtree may be depicted as
$$
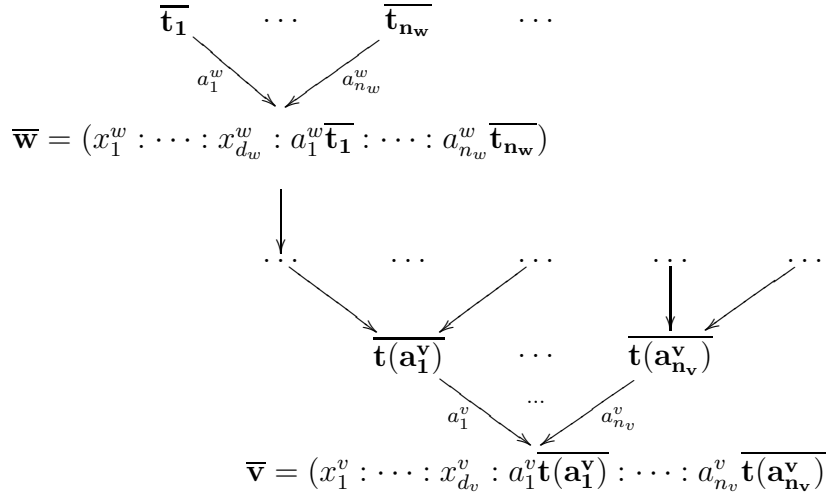
\xymatrix{
\overline{\mathbf{t_1}}\ar[dr]_{a^w_1} & \dots & \overline{\mathbf{t_{n_w}}}\ar[dl]^{a^w_{n_w}} & \dots 
\\
& \save!<0cm,-.4cm>*{
\overline{\mathbf{w}} = (x^w_1:\dots:x^w_{d_w}:a^w_1\overline{\mathbf{t_1}}:\dots:a^w_{n_w}\overline{\mathbf{t_{n_w}}})
}\restore \\
&\ar[d]
\\
& \dots\ar[dr] & \dots & \dots\ar[dl] & \dots\ar[d] & \dots\ar[dl] &
\\
& & \overline{\mathbf{t(a^v_1)}}\ar[dr]_{a^v_1} & \dots\ar@{}[d]|{\dots} & \overline{\mathbf{t(a^v_{n_v})}}\ar[dl]^{a^v_{n_v}} & &
\\
& & & \save!<0cm,-.3cm>*{
\overline{\mathbf{v}} = (x^v_1:\dots:x^v_{d_v}:a^v_1\overline{\mathbf{t(a^v_1)}}:\dots:a^v_{n_v}\overline{\mathbf{t(a^v_{n_v})}}
}\restore& & & 
}$$
This corresponds exactly to a point in $\Graph(v)$, so for any root vertex $v$ in $T$ the restriction of (the orbit of) $V$ to $T(v)$ yields a point in $\Graph(v)$. Now let $v_1$ and $v_2$ be two root vertices with connected rooted subtrees $T({v_1})$ and $T({v_2})$. Assume these rooted subtrees coincide on a subquiver $S$, then this subquiver again is a tree with root vertex $w$. The points in $\mathcal{O}_V$ corresponding to $v_1$ and $v_2$ coincide on all vertices of $S$, thus yielding a point in
$$\Graph(v_1)\times_{\Graph(w)}\Graph(v_2).$$
Repeating this argument until all root vertices are accounted for then precisely yields a point in $\mathcal{F}_n$.
\eop

\textbf{Remark.} These fibers are examples of framed quiver moduli as described by Markus Reineke in \cite{Reineke}. The results here however were obtained independently and through other methods than the results presented in \cite{Reineke}.

\section{Toric Varieties and the Brauer-Severi Fibration}\label{toric}

We can also get a nice description of the fibers using toric geometry. This discussion closely follows \cite{Hille}.

In the set of semistable representations of the quiver we can embed the big torus
$T_Q=(\CC^*)^{\# Q_1}$ as an open subset. On this torus is an action of $T_\alpha=\GL_\alpha =(\CC^*)^{\# Q_0}$, so that the moduli space contains a torus $T=T_Q/T_\alpha=(\CC^*)^{\# Q_1-\# Q_0+1}$ (the extra $1$ comes from the fact that the action is not free).

Let $v = v_0$ be the special source vertex and choose for every other vertex $w$ an arrow $a_w$ from $v$ to $w$. Denote the set of all other arrows by $P$. We can identify $T$ with $(\CC^*)^{\# P}$ by choosing for each point in $T$ the unique representant whose values at the $a_w$ is $1$.

The space in which we are going to construct our fan is then 
$\Hom(\CC^*,T)\otimes \R\cong \R^P$. A well known fact in toric geometry is that one can reconstruct the cones from the variety using one parameter subgroups (1PSG's).

The space of one-parameter subgroups $\Hom(\CC^*,T)=\Z^P$ can be identified with $\Hom(\CC^*,T_Q)/\Image(\Hom(\CC^*,T_\alpha)\to \Hom(\CC^*,T_Q))=\Z^{Q_1}/\Z^{Q_0-1}$.
So every 1PSG $\bar \lambda$ of $T$ corresponds to an equivalence class of vectors $\lambda: Q_1\to \Z$
that differ by a vector of the form
\[
\xi : Q_1 \to \Z: a \mapsto \zeta({h(a)})-  \zeta({t(a)})
\]
where $\zeta: Q_0 \to \Z$ is a character of $T_\alpha$.

\begin{lemma}
In every equivalence class of $\Z^{Q_1}/\Z^{Q_0-1}$ there is a representant with only 
non-negative coefficients.
\end{lemma}
\prf
Suppose that $\lambda \in \Z^{Q_1}$ is a vector with some negative coefficients. 
Choose an arrow $a$ such that $\lambda(a)<0$ and let $V_1$ be the set of vertices that are
targets of paths containing $a$. Let $V_0$ be the complement of this set. Note that $s(a) \in V_0$ because the quiver has no cycles. 
 
Consider the character $\zeta: Q_0 \to \Z$ mapping the vertices in $V_i$ to $i$. This character gives a vector $\xi_a$ that maps every arrow to something non-negative because there are no arrows from $V_1$ to $V_0$. Moreover $\xi_a(a)=1$ so the vector
\[
\lambda' = \lambda - \sum_{a,\lambda(a)<0}\lambda(a)\xi_a 
\]
has no negative entries.
\eop

The previous lemma implies that $\lim_{z\to 0} \bar \lambda(z)$ contains a representation
that assigns to every arrow either a $1$ (for the zero's in the vector) or a $0$ (for the strictly positive values in the vector). Also this  limit representation must be semistable, so there exists a path from $v$ to every vertex $w$ containing only arrows with value $1$. 

Two such different limit points cannot sit in the same $\GL_\alpha$-orbit, because the $\GL_\alpha$-action can never change a zero into a one or vice versa. 

Moreover, every semistable representation with values 0 or 1 can be seen as a limit point of a 1PSG (nl. one coming from a vector with zero's on the arrows with value 1, and positive values on the arrows with value 0). We can conclude that 

\begin{lemma}
There is a one to one correspondence between the limit points of 1PSGs and subquivers of $Q$ which connect $v$ to all other vertices in $Q$.
\end{lemma}
\prf
We identify such subquivers with the representation of $Q$ that maps the arrows of the subquiver to $1$ and the others to $0$. 
\eop

From toric geometry (see e.g. \cite{Fulton}) we know that the poset of cones in the fan of a toric variety 
is isomorphic to the poset of limit points under degeneration. The identification goes as follows: to every limit point $p$ we assign the semigroup
\[
\{\bar\lambda \in \Hom(\CC^*,T):\forall a\in Q_1:p_a=0 \implies (\lim_z\to 0 \bar\lambda(z))_a=0\}
\]
This semigroup comes from the cone
\[
\sigma_p := \{\eta \in \R^P :\forall a\in Q_1:p_a=0\implies
\eta(a)- \sum_{a,\eta(a)<0}\eta(a)\xi_a=0  \}
\]
We know the fiber is smooth so this poset comes from a simplicial set.
Therefore we can conclude 
\begin{theorem}
The set of subquivers of $Q$ such that $v$ can be connected to all other vertices forms a 
poset under the reverse inclusion. This poset is simplicial and the set of cones
\[
\sigma_{Q'} := \{\eta \in \R^P :\forall a\in Q'_1: \eta(a)- \sum_{a,\eta(a)<0}\eta(a)\xi_a=0\}
\]
form a fan.
Its corresponding smooth toric variety is isomorphic to the fiber.
\end{theorem}

We will now use this identifcation to compute the cohomology of the fibers.
For this we use the theorem by Fulton which allows you to compute the cohomology 
ring from the one dimensional cones in the fan:
\begin{theorem}
For a smooth toric variety with fan $\Delta$ the cohomology is given by
the quotient of the polynomial ring generated by the one-dimensional cones $D_i=\N v_i$ with the relations
\begin{itemize}
\item
$D_{i_1}\cdot \dots \cdot D_{i_k}$ if $\N v_{i_1}+\dots +\N v_{i_k}$ is not contained in a cone of $\Delta$.
\item
$\sum_i \langle u,v_i\rangle D_i$ for all possible $u\in \Hom(T,\CC^*)$.
\end{itemize}
\end{theorem}

In order to translate this statement to the setting of this paper, we first of all note that the one-dimensional cones correspond to the subquivers that lack one arrow of the original. However, not all arrows correspond necessarily to a cone because it might be that the representations that map this arrow to zero is not semistable. This happens only if $a=a_w$ for some vertex and no arrow in $P$ terminates in $w$. Let us postpone this case for a moment.

So suppose that there are at least two paths from $v$ to every other vertex. For every arrow $a$ let $D_a$ be the corresponding cone. The vector corresponding to $D_a$ has the following form: if $a \in P$ then $v_a: P \to \Z:b \mapsto \delta_{ab}$, if $a=a_w$ for some vertex $a_w$ 
then
\[
v_a: P \to \Z:b \mapsto \begin{cases}
-1&h(b)=w\\
1&t(b)=w\\
0&\text{else}
\end{cases}
\]

The second set of relations now becomes
\[
\sum_a v_a(b)D_a \text{ for all possible }b\in P.
\] 
These relations imply that for two arrows $a,b$ starting in $v$ and ending in the same vertex $w$, $D_a=D_b$ inside the homology ring. Let us denote this generator by $D_w$. For arrows that do not start in $v$ we have that $D_a = D_{h(a)}-D_{t(a)}$, so the $D_w$ are the generators of the homology ring.
 
Now let us determine what happens with the first set of relations.
A representation is not semistable if there is a vertex that is not the target of a path from
$v$. So the product of some $D_a$ is zero as long a the corresponding set of arrows meets all
paths from $v$ to a certain vertex $w$. As there are arrows from $v$ to every other vertex this also implies that such a set of arrows must meet all the arrows terminating in a given vertex. Therefore the relations can be rewritten as
\[
\prod_{a, h(a)=w} D_a \text{ or } \prod_{a, h(a)=w} (D_{h(a)}- D_{t(a)})  \text{ for all vertices $w$ and $D_v:=0$} 
\]
So we can conclude that

\begin{theorem}
The cohomology ring of the fiber $\pi^{-1}(\xi)$ is isomorphic to the ring
\[
\Z[D_w: w \in Q_0]/(\prod_{a, h(a)=w} (D_{h(a)}- D_{t(a)}):w \in Q_0\setminus \{v\}, D_v)
\]
\end{theorem}
\prf
For the case where there are at least two paths from $v$ to every other vertex 
this follows immediately from the discussion above. 

If there is a vertex $w$ with a unique arrow $a$ that terminates in it, we can construct a new quiver settig by removing this vertex and arrow like this
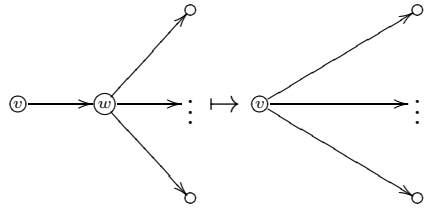
\[
\vcenter{
\xymatrix{
&&\vtx{}\\
\vtx{v}\ar[r]&\vtx{w}\ar[ru]\ar[r]\ar[rd]&\vdots\\
&&\vtx{}}}
\mapsto
\vcenter{\xymatrix{
&&\vtx{}\\
\vtx{v}\ar[rru]\ar[rr]\ar[rrd]&&\vdots\\
&&\vtx{}}}
\]

The modulispace of this new quiver is the same as the old one because the value of $a$
 is invertible and can be set to $1$ using the action in $w$. The homology ring
of the new quiver can be calculated as above, but if one calculates the relations for the old quiver one sees that these relations are the same because $D_w=0$ by the relation in $w$.  
\eop

\end{document}